\documentclass{article}

\usepackage{hyperref}
\usepackage[utf8]{inputenc}
\usepackage[T1]{fontenc}
\usepackage{geometry}
\usepackage{amsmath,amssymb,amsfonts,amsthm}
\usepackage{indentfirst}
\usepackage{graphicx}
\usepackage{verbatim}
\usepackage{mathrsfs}
\usepackage{tikz}
\usetikzlibrary{arrows}
\usepackage[english]{babel}
\usepackage{cite}
\usepackage{natbib}
\usepackage{enumerate}
\usepackage{csquotes}

\newtheorem{thm}{Theorem}[section]
\newtheorem{thm*}{Theorem}

\newtheorem{lem}[thm]{Lemma}

\numberwithin{equation}{section} 

\newcommand{\erdosrenyi}[1]{G_{#1}}
\newcommand{\inhomogeneous}[1]{I_{#1}(k_n)}

\begin{document}

\title{First-order asymptotics for the structure \\ of the inhomogeneous random graph}

%\titlerunning{slowly-varying property}        

\author{Gianmarco Bet \and Kay Bogerd \and Vanessa Jacquier}

\author{
Gianmarco \textsc{Bet}\thanks{University of Florence, Department of Mathematics, Italy. $_{}$\hfill    \href{mailto:gianmarco.bet@unifi.it}
{\texttt{gianmarco.bet@unifi.it}}},\hspace{2pt}
Kay \textsc{Bogerd}\thanks{Utrecht, The Netherlands $_{}$\hfill     \href{mailto:kaybogerd@hotmail.com}
{\texttt{kaybogerd@hotmail.com}} } \hspace{2pt}
Vanessa \textsc{Jacquier}\thanks{ University of Padova, Department of Mathematics, Italy. $_{}$\hfill    \href{mailto:vanessa.jacquier@unipd.it}
{\texttt{vanessa.jacquier@unipd.it}}},\hspace{5pt}
}

%\authorrunning{G.~Bet and V.~Jacquier} 

%\institute{Gianmarco Bet \at University of Florence. Viale Morgagni, 65. 50134 Firenze, Italy \\ Tel.: +39-055-2751491\\ Fax: +39-055-2751452\\ \email{gianmarco.bet@unifi.it}        \\   \and Vanessa Jacquier \at University of Florence. Viale Morgagni, 65. 50134 Firenze, Italy \\ Tel.: +39-055-2751467\\ Fax: +39-055-2751452\\ \email{vanessa.jacquier@unifi.it} }

%\date{Received: date / Accepted: date}
% The correct dates will be entered by the editor

\maketitle

\begin{abstract}
In the inhomogeneous random graph model, each vertex $i\in\{1,\ldots,n\}$ is assigned a weight $W_i\sim\text{Unif}(0,1)$, and an edge between any two vertices $i,j$ is present with probability $k(W_i,W_j)/\lambda_n\in[0,1]$, where $k$ is a positive, symmetric function and $\lambda_n$ is a scaling parameter that controls the graph density. When $\lambda_n=1$ (resp.~$\lambda_n=\Theta(n)$) the typical resulting graph is dense (resp.~sparse). The goal of this paper is the study of structural properties of \textit{large} inhomogeneous random graphs. We focus our attention on graph functions that grow sufficiently slowly as the graph size increases. Under some additional technical assumptions, we show that the first-order asymptotic behavior of all such properties is the same for the inhomogeneous random graph and for the ER random graph. Our proof relies on two couplings between the inhomogeneous random graph and appropriately constructed ER random graphs. We demonstrate our method by obtaining asymptotics for two structural properties of the inhomogeneous random graph which were previously unknown. In the sparse regime, we find the leading-order term for the maximum degree. In the dense regime, we find the asymptotics of the $\gamma$-quasi-clique number. 

\vspace{0.1cm}
\noindent \emph{Keywords}: Inhomogeneous random graphs, Erd\H{o}s and R\'enyi graphs; coupling methods, maximum degree, quasi-cliques, sparse and dense regimes, graph asymptotics.

% \PACS{PACS code1 \and PACS code2 \and more}
\vspace{0.1cm}
\noindent {\emph{MSC codes:} 05C80, 05C35, 60C05, 60F05.}
\end{abstract}

\section{Introduction}

\noindent In the last decade, random graph models have emerged as a powerful tool to understand natural and artificial  networks alike. Random graphs are now routinely used to model social, biological and communication networks, among others. Additionally, random graph models present profound theoretical challenges, as evidenced by numerous recent developments, see for example the overviews \cite{chatterjee2017large,goldschmidt2017scaling,lovasz2012large}.
Perhaps the very first random graph model dates back to the 1950s and the 1960s, when Erd\H{o}s and R\'enyi \cite{erdHos1960evolution}, and independently Gilbert \cite{gilbert1959random}, considered a model in which all graphs with $n$ vertices and $m$ edges are assigned equal probability. This is now known as \textit{Erd\H{o}s-R\'enyi} (ER) \textit{random graph}, in recognition of the key role that Paul Erd\H{o}s and Alfr\'ed R\'enyi played in understanding the typical structure of the resulting graph as $m$ varies. The ER random graph is more commonly formulated as follows: take $n$ vertices and, for any two vertices, draw an edge between them with some fixed probability $p\in(0,1)$. Hence the degree distribution of each vertex in the ER random graph is the same, and it follows a binomial distribution $\text{Bin}(n,p)$. In particular, the degree of any vertex is highly concentrated around its mean. These features make the ER random graph a \textit{homogeneous} model. On the other hand, real-world complex networks are highly inhomogeneous. For example, it is observed that the empirical degree distribution of real-world networks often decays like a power-law. More generally, the homogeneous nature of the ER random graph makes it unsuitable to model real-world complex networks. Because of this, during the last decades a great effort has been devoted to constructing generalizations of the ER random graph that reproduce the properties that are more often observed in real-world networks. Examples include the configuration model \cite{dhara2017critical}, the preferential attachment model \cite{deijfen2009preferential,baldassarri2021asymptotic}, the inhomogeneous random graph \cite{bollobas2007phase} and the rank-1 random graph \cite{bet2020big,bhamidi2010scaling,bhamidi2012novel}, as well as several versions of the geometric random graph (GRG) \cite{krioukov2016clustering,bet2021detecting,penrose2003random}. See also \cite{van2016random} for an overview.

In this paper we focus on the inhomogeneous random graph (IRG). In this model, there are $n$ vertices $V=[n]:=\{1,\ldots, n\}$, and each vertex $i$ is assigned a uniform independent weight $W_i\sim\text{Unif}(0,1)$. The connectivity structure is described by a symmetric, positive function $k:[0,1]^2\to[0,1]$, known as the \textit{kernel} and a sequence of scaling parameters $(\lambda_i)_{i=1}^n$. Any two vertices $i\neq j\in[n]$ are then connected with probability $k(W_i,W_j)/\lambda_n$. Note that usually in the literature the connection probability are defined as $\min \{k(W_i,W_j)/\lambda_n, 1\}$. In our setting $k$ is a bounded function, so without loss of generality we assume that $k(W_i,W_j)/\lambda_n$ is less than one.
The resulting topological structure is rather complicated in general and depends crucially on the precise form of $k$, see for example \cite{bollobas2007phase}. Moreover, when $\lambda_n \equiv 1$, the typical resulting graph is dense. On the other hand, if $\lambda_n\to\infty$, then the edge density depends on the precise growth rate of $\lambda_n$, for example when $\lambda_n = O(n)$, the resulting typical graph has $O(n)$ edges. We investigate the structural properties of the IRG as the number of vertices $n$ grows. Instead of focusing on a specific graph property as is often done \cite{devroye2014connectivity,mckinley2019superlogarithmic,shang2013large}, we derive general asymptotic formulas for a class of graph properties in terms of the corresponding asymptotics for the ER random graph. 

The broad idea of studying certain classes of random graphs by building a two-sided coupling with the ER random graph has recently been developed in \cite{achlioptas2018symmetric,bangachev2025sandwiching,liu2022testing}.  In \cite{liu2022testing}, the coupling is used to prove that when the underlying GRG dimension is sufficiently large, the sparse ER random graph \textit{cannot} be statistically distinguished from a GRG. The works \cite{achlioptas2018symmetric,bangachev2025sandwiching} are perhaps closer to ours. In \cite{bangachev2025sandwiching} the authors construct an explicit coupling between a GRG and two ER random graphs in order to study, among other things, the critical thresholds of edge-monotone properties of the GRG. Finally, in \cite{achlioptas2018symmetric} the authors give very general sufficient condition for the existence of a coupling between a \textit{uniform} measure on an arbitrary subset of the set of all graphs with $n$ vertices and two ER random graphs.

Roughly speaking, we focus only on those graph properties that are determined by the local structure of the graph. For example, the maximum degree $\Delta$ and the $\gamma$-quasi clique satisfy our assumptions, but the average distance does not. 
In particular, we will apply our method to prove that the asymptotic behavior of the $\gamma$-quasi clique and the maximum degree for an inhomogeneous random graph is the same of an ER random graph. 
For the precise definition of these graph properties, see sections \ref{sec:Mdegree} and \ref{sec:quasi-clique}. 
In the following theorems we assume that $(\inhomogeneous{n})_{n=1}^{\infty}$ is a sequence of inhomogeneous random graphs with kernel $k_n = k/\lambda_n$ and there exists $(m,m)\in[0,1]^2$ such that $k(x,y)\leq k(m,m)$ for all $x,y\in[0,1]$. 
In particular $k$ is a bounded function, which excludes some interesting kernels such as scale-free IRGs.
Further assume that $k$ is continuous \textit{at} $(m,m)$ and set $p^{\text{max}}_n := k(m,m)/\lambda_n$.

\begin{thm}[Maximum degree $\Delta$ of sparse inhomogeneous random graphs] \label{thm:degree}
Under the assumptions above, if $\lambda_n= cn$ with $c\in\mathbb{R}^+$ such that $c>k(m,m)$, then for any $\varepsilon>0$ we have 
\begin{align}%
\mathbb P \left (\Delta(\inhomogeneous{n}) \in [(1-\varepsilon)\Delta_{n}, (1 + \varepsilon)\Delta_{n}] \right ) \to 1,
\end{align}%
as $n\to\infty$, where 
$$\Delta_{n}:=\log (n)/\log(\log(n)).$$
\end{thm}

\begin{thm}[$\gamma$-quasi clique of dense inhomogeneous random graphs]\label{thm:cliques}
Under the assumptions above, if $\lambda_n \equiv 1$ and $p^{\text{max}}_n <\gamma<1$, then for every $\varepsilon>0$ we have
\begin{equation}
P(\omega^\gamma(\inhomogeneous{n}) \in [(1-\varepsilon)\omega_n^\gamma, (1 + \varepsilon)\omega_n^\gamma]) \to 1,
\end{equation} 
as $n \to \infty$, where $\omega_n^\gamma:=2\log(n)/ D(\gamma,p^{\text{max}}_n)$ and $D(\gamma,p^{\text{max}}_n)$ is the Kullback-Leibler divergence between the Bernoulli distributions \emph{Ber}$(\gamma)$ and \emph{Ber}$(p^{\text{max}}_n)$.
\end{thm}%

We describe the asymptotics for the IRG in terms of two ER random graphs by using a coupling argument: a global upper bound, and a local lower bound. The latter is constructed such that only those vertices with weights close to the argmax of the kernel are present. Our result formalizes the intuitive idea that some structural properties of the IRG are characterized by the densest regions of the graph. We give precise conditions under which this happens. Crucially, our technique works both in the dense case $\lambda_n\equiv 1$ and in the sparse case $\lambda_n=O(n)$. To illustrate this, we study the maximum degree in the sparse IRG, and the quasi-clique number in the dense IRG. However, our technique is only able to capture the leading order behavior of such properties. Obtaining higher-order asymptotics would require a more careful analysis of the behavior of $k$ in the whole domain $[0,1]^2$, not only close to its argmax.

The structure of the ER random graph for large $n$ and all values of $p\in[0,1]$ is well understood owing to a rich literature on the subject. Recent results presented in \cite{de2022unusually} describe the size of unusually large components in near-critical ER graphs, providing a useful benchmark for analyzing analogous structural properties in inhomogeneous random graphs. 
On the contrary, few results on the structure of the IRG are available. 
In particular, we refer to \cite{bollobas1998random, frieze2015introduction, frieze2024maximum} for the study of the maximum degree on the ER random graph.
Finally, \cite{balister2019dense, matula1976largest, matula1972employee} studies the size of the largest clique and quasi-clique in an ER random graph. In \cite[Section 3.1]{bogerd2020cliques} a heuristic is given to extend these results to the IRG. In fact, our work was inspired by the recognition that the heuristic in \cite{bogerd2020cliques} is rather general and can be turned into a general technique to study graph properties other than the clique number. 

The rest of the paper is organized as follows. In Section \ref{sec:model_main_results} we introduce the IRG formally and we give a simple lemma, and in Section \ref{sec:coupling_IRG_ER} we apply it to obtain our main result. This gives the asymptotics for properties of the IRG in terms of the corresponding properties of the ER random graph. In Section \ref{sec:applications} we apply our general result to obtain the first-order asymptotic behavior of various graph properties. For the sparse IRG, we study the maximum degree in Section \ref{sec:Mdegree}. For the dense IRG, we study the quasi-clique number in Section \ref{sec:quasi-clique}. Finally, in Section \ref{sec:discussion}, we give an example of a property of the IRG that we are not able to describe in terms of the corresponding ER property, and we discuss an issue arising when trying to generalize our result.

\paragraph{Notation.}
In this paragraph, we introduce the notation used in the rest of this paper. 
We denote by \emph{log} the natural logarithm. 
Given a random variable $U$ with distribution $\mathcal{D}$, we write $ U \sim \mathcal{D}$. 
Moreover, given two functions $f, \, g : A \subseteq \mathbb{R} \to \mathbb{R}$, we use standard asymptotic notation: $f_x=O(g_x)$ when $f_x/g_x$ is bounded for $x \to \infty$, and 
$f_x \gg g_x$ to mean $g_x=o(f_x)$. 
Furthermore, if $A=\mathbb{N}$, we write $f_n \asymp g_n$ to indicate that $\lim_{n\to\infty}f_n/g_n=1$.
Finally, we say that a sequence of events $(\mathcal A_n)_{n=1}^{\infty}$ holds \textit{with high probability} if $\lim_{n\to\infty}\mathbb P(\mathcal A_n) = 1$.

\section{Model description}\label{sec:model_main_results}
A graph $G=(V,E)$ is a finite set of vertices $V=[n]:=\{1,\ldots,n\}$ and edges $E\subseteq V\times V$ between them. A graph is \textit{undirected} if for all $(i,j)\in E$, also $(j,i)\in E$. From now on we will only consider undirected graphs unless specified otherwise. Let $\mathcal{G}_n$ denote the set of graphs with vertex set $V=[n]$, and let $\mathcal G := \bigcup_{n=1}^{\infty}\mathcal G_n$. A \emph{kernel} $k$ is a  positive, symmetric, measurable function $k: [0,1]^2 \mapsto [0,1]$ such that $k(x,y) = k(y,x)$ for all $x,y\in[0,1]$. The inhomogeneous random graph $\inhomogeneous{n}$ with vertex set $V:=[n]$ and with kernel  $k_n:=k/\lambda_n$ is constructed as follows. Each vertex $i \in V$ is assigned a \emph{weight} $W_i \sim \text{Unif}(0,1)$. Any two vertices $i,j\in V$ are connected with probability
\begin{align}\label{eq:connection_probability_IRG}%
p_{ij}^{(n)} := \mathbb P((i,j)\in E\vert (W_i)_{i=1}^n) = \frac{k(W_i,W_j)}{\lambda_n},
\end{align}%
where $(\lambda_n)_{n=1}^{\infty}$ is a sequence of positive numbers. Note that when $\lambda_n\equiv C>0$, $p_{ij}^{(n)}$ is independent of $n$.  In particular, the resulting graph is dense, i.e., the average number of total edges grows as $O(n^2)$ as the number of vertices $n$ tends to infinity. On the other hand, when, say, $\lambda_n = C n$ for $C>0$, the average number of edges in the graph grows as $O(n)$, i.e., the resulting graph is sparse. We have tacitly assumed that $k_n$ is a bounded function. 
One may assume that $k_n$ is not symmetric, in which case the resulting random graph is \textit{directed}. Our technique also works in this more general setting at the cost of additional technical assumptions, but for simplicity we work in the simpler setting defined above. The sparse inhomogeneous random graph is a generalization of several well-known random graph models. When $k(x,y)\equiv p$ for all $x,y\in[0,1]$, the resulting IRG coincides with the ER random graph on $n$ vertices and with connection probability $p_n=p/\lambda_n$, which we denote as $\erdosrenyi{n,p_n}$. When $k(x,y) = k_1(x)k_2(y)$, the resulting IRG is the rank-1 model  \cite{bollobas2007phase, van2013critical, bhamidi2010scaling, bhamidi2012novel, bet2020big}. Finally, when $k(x,y)\in\mathcal S$ for all $x,y\in[0,1]$ and $\mathcal S\subset [0,1]$ is a finite set, we obtain the well-known \textit{stochastic block model} \cite{decelle2011asymptotic, rohe2011spectral}.

A subgraph $G'$ of $G\in\mathcal G_n$ is a graph with vertex set $ V' \subseteq V$ and edge set $E' \subseteq E$, and we indicate it with $ G' \subseteq G$. When $V'= V$ and $G' \subseteq G$, we denote it as $G'\preccurlyeq G$. Note that $\preccurlyeq$ defines a partial order relation on $\mathcal G_n$. We say that a function $f:\mathcal G\to \mathbb R^+$ is a \textit{property} of the graphs in $\mathcal G$. For example, if $G=(V,E)$ and $f(G) := \vert E\vert/\binom{n}{2}$, then $f$ represents the density of edges. We are interested in deriving asymptotics for $f(\inhomogeneous{n})$ as $n\to\infty$, when $f$ belongs to a class of properties to be specified. Our main result hinges on the following idea. 
\begin{lem}\label{lem:general_coupling_lemma}
For every $n\geq1$, let $G_n, G_n', G_n''$ be $\mathcal{G}_n$-valued random variables defined on the same probability space $(\Omega,\mathcal F, \mathbb P)$ and such that 
\begin{align}\label{eq:coupling_condition}%
G'_n \preccurlyeq G_n \preccurlyeq G''_n, 
\end{align}%
almost surely for all $n\geq1$. Let $f: \mathcal G \rightarrow \mathbb R$ be a graph property which is non-decreasing with respect to the partial order relation $\preccurlyeq$ defined above. Further, let $f_n: \mathbb{N}  \rightarrow \mathbb{R}$ be a function such that, for any small $\varepsilon>0$,
\begin{align}\label{eq:upper_lower_bounds_condition}%
\mathbb{P}(f(G'_n) > (1-\varepsilon) f_n) \to 1, \quad \mathbb{P}(f(G''_n) < (1+\varepsilon) f_n) \to 1,
\end{align}
as $n \to \infty$. Then, for any $\varepsilon>0$,
\begin{equation}\label{asymptoticbehavior}
    \mathbb{P}(f(G_n) \in ((1-\varepsilon) f_n, (1+\varepsilon) f_n)) \to 1, 
\end{equation}
as $n\to\infty$.
\end{lem}
\begin{proof}%
Since $G'_n \preccurlyeq G_n \preccurlyeq G''_n$ for all $n\geq1$, we have $f(G'_n) \leq f(G_n) \leq f(G''_n)$ almost surely, because $f$ is non-decreasing. By \eqref{eq:upper_lower_bounds_condition} we have
\begin{equation}
  (1-\varepsilon) f_n < f(G'_n) \leq f(G_n) \leq f(G''_n) < (1+\varepsilon) f_n,
\end{equation}
with high probability for any $\varepsilon>0$.
\end{proof}%
In Lemma \ref{lem:general_coupling_lemma}, one should interpret $f(G_n)$ as the quantity of interest, and $G_n'$, $G_n''$ as simpler models for which the asymptotic behavior of $f(G_n')$, $f(G_n'')$ is already known.
Lemma \ref{lem:general_coupling_lemma} is rather general and, indeed, also rather straightforward. However, difficulties arise when applying it to concrete settings. First, one must couple $G_n$ to two simpler models $G_n'$ and $G_n''$ such that \eqref{eq:coupling_condition} is satisfied almost surely. Second, one must obtain the asymptotic lower bound (resp.~upper bound) \eqref{eq:upper_lower_bounds_condition} for the random graph $G_n'$ (resp.~$G_n''$) and make sure that the two bounds match, that is, the two limits in \eqref{eq:upper_lower_bounds_condition} hold with the same $f_n$.

In the rest of the paper, we will take $G_n:=\inhomogeneous{n}$ to be an inhomogeneous random graph, and $f$ some property that we are interested in. We will focus on constructing a coupling between $G_n$ and two appropriate ER random graphs $G_n':=\erdosrenyi{n,p'}$, $G_n'':=\erdosrenyi{n,p''}$ such that \eqref{eq:coupling_condition} is satisfied. Note that the two probabilities $p'$ and $p''$ may depend on $n$, but we omit it to avoid notational clutter. We will then rely on results from the literature on the asymptotics of properties of the ER random graph to verify that \eqref{eq:upper_lower_bounds_condition} holds and thus apply Lemma \ref{lem:general_coupling_lemma}. This is not straightforward because the asymptotic behavior of $f(G_n')$ and $f(G_n'')$ in general depends on the connection probabilities $p'$ and $p''$ as well as the graph size $n$. If this is the case, then we have $f(G_n')\asymp f_{n,p'}$ and $f(G_n'')\asymp f_{n,p''}$, and so the upper and lower bounds do not match unless $p'$ is close to $p''$. Hence, some of our work is dedicated to finding some $\tilde p\in[0,1]$ such that $f(G_n')\asymp f_{n,\tilde p} \asymp f(G_n'')$. From this will follow that \eqref{eq:upper_lower_bounds_condition} holds with $f_n:=f_{n,\tilde p}$.

Finally, we note that Lemma \ref{lem:general_coupling_lemma}, as well as the rest of our results below, also hold if the graph property $f$ is non-increasing with respect to $\preccurlyeq$, as long as the remaining assumptions are modified accordingly. However, we were not able to find non-trivial graph properties that satisfy the resulting assumptions, and thus we did not pursue this further.

\subsection{Coupling the inhomogeneous random graph and the Erd\H{o}s-R\'enyi random graph}\label{sec:coupling_IRG_ER}
We are now ready to investigate the asymptotic behavior of $f(\inhomogeneous{n})$. In Section \ref{subsec:assumptions} we describe in detail all our assumptions on $f$ and we collect our results in a theorem, which we then prove in Section \ref{subsec:coupling}. 
\subsubsection{Assumptions and main result}\label{subsec:assumptions}
Typically, the asymptotic behavior of a property $f(\erdosrenyi{n,p^*})$ of an ER random graph $\erdosrenyi{n,p^*}$ depends on the number of vertices $n$, as well as the connection probability $p^*$. Note that we allow $p^*$ to depend on $n$. Formally, we assume that the following holds:
\begin{enumerate}[({H}1)]
\item If $G_1$, $G_2$ are two graphs such that $G_1 \preccurlyeq G_2$, then $f(G_1) \leq f(G_2)$. In other words, $f$ is non-decreasing with respect to the partial order relation $\preccurlyeq$. \label{item:f_increasing_partial_order}
\item Let $G_1=([n],E)$ and $G_2=([n+k],E)$ be two graphs such that $G_2$ is obtained from $G_1$ by adding $k$ zero-degree vertices. Then $f(G_1) \geq f(G_2)$. In other words, $f$ is non-increasing with respect to the number of vertices with degree zero. \label{item:f_decrasing_zero_degree}
\item The asymptotic behavior of $f(G_{n,p^*})$ as $n\to\infty$ is known. Formally, let $(\erdosrenyi{n,p^*})_{n=1}^{\infty}$ be a sequence of ER random graphs such that $\erdosrenyi{n,p^*}$ has $n$ vertices and connection probabiltiy $p^*=p^*_n$. We assume that there exists some $f_{n,p}: \mathbb{N} \times[0,1]\rightarrow \mathbb{R}$ such that for any $\varepsilon>0$,
\begin{align}\label{eq:asymptotic_property_Erdos_Renyi}%
\mathbb{P}((1-\varepsilon) f_{n,p^*} < f(\erdosrenyi{n,p^*}) < (1+\varepsilon) f_{n,p^*}) \to 1     
\end{align}%
as $n\to\infty$.  \label{item:asymptotic_property_Erdos_Renyi}
\end{enumerate}
Assumption (H\ref{item:f_decrasing_zero_degree}) is technical, but is crucial for our proofs.
Many results that take the form of equation \eqref{eq:asymptotic_property_Erdos_Renyi} are available in the literature in \cite{ frieze2015introduction, frieze2024maximum} for the maximum degree and in \cite{ bollobas1976cliques, bogerd2020cliques, dolevzal2017cliques, matula1972employee, matula1976largest} for the cliques. When $p\in(0,1)$, we assume that $f_{n,p}$ satisfies the following properties:
\begin{enumerate}[({A}1)]
\item The function $p \mapsto f_{n,p}$ is continuous, uniformly in $n \in \mathbb{N}$, i.e.,
\begin{equation}
    \forall p_0, \, \forall \varepsilon>0, \, \exists \, U_{\varepsilon}(p_0): \,
    \forall p \in U_{\varepsilon}(p_0), \,\,  |f_{n,p}-f_{n,p_0}|<\varepsilon,\quad  \forall n\in\mathbb N.
\end{equation} \label{item:property_continuous}
\vspace{-2em}
\item The function $n\to f_{n,p}$ is non-decreasing %increasing, 
for each $p \in (0,1)$. \label{item:property_increasing}
\item There exists a function $r: \mathbb{N} \rightarrow \mathbb{R}$ such that $nr(n) \to \infty$, $r(n)\to0$ as $n \to \infty$ and
\begin{equation}\label{limfunctions}%
\lim_{n \to \infty}\sup_{p\in(0,1/\lambda_n)} \left |\frac{f_{ n r(n),p}}{f_{n,p}}-1 \right |=0,
\end{equation}
where $\lambda_n$ is the scaling factor in the inhomogeneous random graph kernel $k_n=k/\lambda_n$.
\label{item:property_slowly_varying}
\end{enumerate}
% In the sparse regime $p=c/\lambda_n$, the function $f_{n,p}$ depends on both $n$ and $c$, so we denote it by $f_{n,c}$. We then assume that the function $(n,c) \mapsto f_{n,c}$ satisfies assumptions analogous to those introduced above. In particular, (A\ref{item:asymptotic_property_Erdos_Renyi}) can be restated in this setting as
% %
% \begin{equation}\label{limfunctions_sparse}%
% \lim_{n \to \infty} \sup_{c \in \mathbb{R}} \left |\frac{f_{ n r(n),c}}{f_{n,c}}-1 \right |=0.
% \end{equation}
%

Assumptions (A\ref{item:property_continuous}), (A\ref{item:property_increasing}) are satisfied by most concrete examples and thus, in effect, are not restrictive. On the other hand, assumption (A\ref{item:property_slowly_varying}) is not always satisfied but is crucial for our proofs, see Section \ref{sec:discussion}. Condition (A\ref{item:property_slowly_varying}) implies that $f(G_{n,p})$ grows sufficiently slowly so that $f(G_{n,p})\asymp f(G_{n',p})$ for a (significantly) smaller graph $G_{n',p}\subseteq G_{n,p}$. Intuitively, this means that $f$ is a property that only depends on the local structure of the graph. For example, if $f_{n,p_n} = \log(n)^kg(p_n)$, for $k\geq1$ and $g(p_n)$ a function of $p_n$ with $p_n=p$ ($\lambda_n=1$), then (A\ref{item:property_slowly_varying}) is satisfied by taking $r(n) = 1/\log({n})^{\alpha}$ for any $\alpha \geq 1$. 

We are finally able to state our main result describing the asymptotics for the structure of the IRG in terms of the corresponding asymptotics for the ER random graph.
\begin{thm}\label{thm:asymptotics_IRG}%
Let $(\inhomogeneous{n})_{n=1}^{\infty}$ be a sequence of inhomogeneous random graphs with kernel $k_n = k/\lambda_n$. Assume that there exists $(m,m)\in[0,1]^2$ such that $k(x,y)\leq k(m,m)$ for all $x,y\in[0,1]$.  Further assume that $k$ is continuous \textit{at} $(m,m)$, and set $p^{\text{max}}_n := k(m,m)/\lambda_n$. Let $f:\mathcal G\to\mathbb R^+$ be a graph property that satisfies (H\ref{item:f_increasing_partial_order})--(H\ref{item:asymptotic_property_Erdos_Renyi}) and (A\ref{item:property_continuous})--(A\ref{item:property_slowly_varying}). Then, for any $\varepsilon>0$,
\begin{align}%
\mathbb P((1-\varepsilon)f_{n,p^{\text{max}}_n} \leq f(\inhomogeneous{n})\leq (1+\varepsilon) f_{n,p^{\text{max}}_n} ) \to 1,
\end{align}%
as $n\to\infty$.
\end{thm}%
The assumption that $k$ attains its maximum value in a point of the diagonal $(m,m)\in[0,1]^2$ (instead of, say, in $(\bar x, \bar y)\in[0,1]^2$ with $\bar x\neq \bar y$) is a technical assumption that allows us to exploit known results on ER random graphs, see the discussion in Section \ref{sec:discussion}.

\subsubsection{The couplings}\label{subsec:coupling}
In this section we prove Theorem \ref{thm:asymptotics_IRG} by leaning on Lemma \ref{lem:general_coupling_lemma}. To this end, we couple $\inhomogeneous{n}$ to two ER random graphs $\erdosrenyi{n,p'}$, $\erdosrenyi{n,p''}$ with appropriate connection probabilities such that \eqref{eq:coupling_condition} and \eqref{eq:upper_lower_bounds_condition} are satisfied. We first construct the upper bound $\erdosrenyi{n,p''}$ and then the lower bound $\erdosrenyi{n,p'}$, since the latter is quite harder than the former. 
After this, we will show how to apply Theorem \ref{thm:asymptotics_IRG} to several concrete examples in Section \ref{sec:applications}. 
\paragraph{Upper bound.}  From now on, we tacitly assume that all random variables are defined on the same probability space $(\Omega,\mathcal F,\mathbb P)$. Let $(U_{ij})_{i,j=1}^n$ and $(W_i)_{i=1}^n$ be families of independent random variables such that $U_{ij}\sim\text{Unif}(0,1)$ and $W_i\sim\text{Unif}(0,1)$. Let $V:=[n]$, and
\begin{align}\label{coupling0}
E &:=\{(i,j) \in V \times V \, | \, U_{ij} \leq k(W_i,W_j)/\lambda_n\}, \notag \\
E''&:=\{(i,j) \in V \times V \, | \, U_{ij} \leq p^{\text{max}}_n\}.
\end{align}
It is straightforward to check that $G_n:=(V, E)$ is distributed as the inhomogeneous random graph $\inhomogeneous{n}$, and $G_n'':=(V,E'')$ is distributed as the ER random graph $\erdosrenyi{n, p^{\text{max}}_n}$. Since $k(W_i,W_j)/\lambda_n \leq p^{\text{max}}_n$ almost surely for every $i\neq j \in [n]$, it follows that $G_n\preccurlyeq G_n''$ almost surely. Hence $G_n$ and $G_n''$ verify the right-most inequality in \eqref{eq:coupling_condition}. 

By (H\ref{item:f_increasing_partial_order}) we have $f(G_n) \leq f(G_n'')$ almost surely. Furthermore, by (H\ref{item:asymptotic_property_Erdos_Renyi}),
\begin{equation}\label{eq:asymptotic_upper_bound}%
f(G_n) \leq f(G_n'') \leq (1+ \varepsilon) f_{n,p^{\text{max}}_n},
\end{equation}%
with high probability for any fixed $\varepsilon>0$. This shows that $\mathbb{P}(f(G_n) \leq (1+\varepsilon) f_{n,p^{\text{max}}_n}) \to 1$ as $n\to\infty$, so the upper bound condition in \eqref{eq:upper_lower_bounds_condition} is satisfied with $f_n:= f_{n,p^\text{max}_n}$.

\paragraph{Lower bound.} Next we construct the coupling that provides the lower bound in \eqref{eq:coupling_condition} and \eqref{eq:upper_lower_bounds_condition}. Take the same families of random variables $(U_{ij})_{i,j=1}^n$ and $(W_i)_{i=1}^n$ as above. Again we let $V:=[n]$, $E :=\{(i,j) \in V \times V \, | \, U_{ij} \leq k(W_i,W_j)/\lambda_n\}$ and $G_n:=(V,E)$, so that $G_n$ is distributed as $\inhomogeneous{n}$.

Let $n\mapsto r(n)$ be the function given by (A\ref{item:property_slowly_varying}). Set
\begin{align}%
R_n := [m-r(n),m+r(n)] \cap [0,1],
\end{align}%
and
\begin{align}%
S_n := \{(i,j) \in V \times V~\vert~ (W_i,W_j) \in R_n^2\}.
\end{align}%
In other words, $S_n$ is the subset of (possible) edges such that both ends have weight close to $m$. Note that $R_n= [m-r(n),m+r(n)]$ when $n$ is large. Finally, let $G'_n:=(V,E')$, with
\begin{align}\label{coupling00}
E'=\{(i,j) \in S_n~\vert~ U_{ij} \leq p^{\text{inf}}_n \},
\end{align}%
where $p^{\text{inf}}_n := \inf_{(x,y)\in R_n^2} k(x,y)/\lambda_n$. 
By construction, $p^{\text{inf}}_n \leq k(W_i,W_j)/\lambda_n$ whenever $(W_i,W_j)\in R_n^2$. Since the only possible edges in $G_n'$ are those in $S_n$, it follows that $G_n'\preccurlyeq \inhomogeneous{n}$ almost surely. Hence the left-hand side of \eqref{eq:coupling_condition} holds. 

Next we verify the left-hand side of \eqref{eq:upper_lower_bounds_condition}. The crucial issue is that $G_n'$ is \textit{not} distributed as an ER random graph with connection probability $p^{\text{inf}}_n$, because all the edges that do not lie in $S_n$ have zero probability of being present in $G_n'$. To circumvent this, we construct an ER random graph $\tilde G_n := (V', E')$, with
\begin{align}\label{eq:vertex_set_graph_lower_bound}%
V' = \{i\in V ~\vert~ W_i \in R_n \}.
\end{align}%
Conditionally on $(W_i)_{i=1}^n$, $\tilde G_n$ is distributed as the ER random graph $\erdosrenyi{\vert V'\vert,p^{\text{inf}}_n}$, where $\vert V'\vert$ denotes the cardinality of $V'$. By (H\ref{item:f_decrasing_zero_degree}), we have $f(G_n') \geq f(\tilde G_n)$  almost surely, so we are left to study the asymptotic behavior of $f(\tilde G_n)$. By (H\ref{item:asymptotic_property_Erdos_Renyi}), for any $\varepsilon>0$, $\mathbb{P}((1-\varepsilon) f_{\vert V'\vert,p^{\text{inf}}_n} < f(\tilde G_n)~\vert~(W_i)_{i=1}^n) \to 1$ almost surely as $n\to\infty$. By the dominated convergence theorem,
\begin{align}\label{eq:lower_bound_ER}%
\mathbb E[\mathbb{P}((1-\varepsilon) f_{\vert V'\vert,p^{\text{inf}}_n} < f(\tilde G_n)~\vert~(W_i)_{i=1}^n)] = \mathbb{P}((1-\varepsilon) f_{\vert V'\vert,p^{\text{inf}}_n} < f(\tilde G_n)) \to1,
\end{align}%
as $n\to\infty$. Next we lower bound $f_{\vert V'\vert,p^{\text{inf}}_n}$ to obtain an asymptotic term that matches the upper bound \eqref{eq:asymptotic_upper_bound}. Without loss of generality, choose $\varepsilon>0$ such that $3\varepsilon < 1$. Since $k$ is continuous in $(m,m)\in[0,1]^2$ 
and by (A\ref{item:property_continuous}), we have
\begin{align}\label{dis03}%
(1-\varepsilon)f_{\vert V'\vert,p^{\text{inf}}_n} \geq (1-2\varepsilon)f_{\vert V'\vert,p^{\text{max}}_n}
\end{align}
when $n$ is large. To circumvent the issue that $V'$ is random, we exploit the fact that it is highly concentrated. By Hoeffding's inequality \cite[Theorem 2.8]{boucheron2013concentration}, for any $t\geq0$ we have
\begin{equation}\label{Hineq0}
\mathbb{P}(\vert V'\vert \leq \mathbb{E}[\vert V'\vert]-t) \leq \exp(-2t^2/n),
\end{equation}
where the expected value of $\vert V'\vert$ is 
\begin{align}\label{ESn}
\mathbb{E}[\vert V'\vert] = n \mathbb{P}(W_i \in R_n) =2nr(n).
\end{align}%
By (A\ref{item:property_increasing}) it follows that, with high probability,
\begin{align}%
(1-2\varepsilon)f_{\vert V'\vert,p^{\text{max}}_n} \geq (1-2\varepsilon)f_{\mathbb E[\vert V' \vert] - t,p^{\text{max}}_n} = (1-2\varepsilon)f_{2nr(n) - t,p^{\text{max}}_n} \geq (1-2\varepsilon)f_{nr(n),p^{\text{max}}_n},
\end{align}%
since $2nr(n) - t \geq nr(n)$ when $n$ is large. Using (A\ref{item:property_slowly_varying}) we obtain
\begin{align}\label{dis00}%
(1-2\varepsilon)f_{nr(n),p^{\text{max}}_n} \geq (1-3\varepsilon)f_{n,p^{\text{max}}_n}.
\end{align}
From \eqref{dis03}--\eqref{dis00} follows that
\begin{align}\label{dis0}%
(1-\varepsilon)f_{\vert V'\vert,p^{\text{inf}}_n} \geq (1-3\varepsilon)f_{n,p^{\text{max}}_n},
%(1-\varepsilon)f_{n,p^{\text{max}}_n} \leq \Big(1-\frac{\varepsilon}{4}\Big)f_{|V'|,p^{\text{inf}}_n}.
\end{align}
with high probability as $n\to\infty$. Putting toghether \eqref{eq:lower_bound_ER} and \eqref{dis0} we get
\begin{align}%
f(G_n')  \geq f(\tilde G_n) \geq (1- \varepsilon) f_{\vert V'\vert,p^{\text{inf}}_n} \geq (1-3\varepsilon) f_{n,p^{\text{max}}_n}
\end{align}%
with high probability as $n\to\infty$.
Since $\varepsilon <1/3$ was chosen arbitrarily, the left-hand side of \eqref{eq:upper_lower_bounds_condition} is satisfied with $f_n := f_{n,p^{\text{max}}_n}$.

We have verified that all the conditions of Lemma \ref{lem:general_coupling_lemma} hold, and this concludes the proof of Theorem \ref{thm:asymptotics_IRG}.

\section{Applications for dense regime}\label{sec:applications}
In this section, we consider two concrete examples of properties $f$ that satisfy (H\ref{item:f_increasing_partial_order}), (H\ref{item:f_decrasing_zero_degree}) and whose asymptotic behavior for the dense ER random graph is known, i.e., (H\ref{item:asymptotic_property_Erdos_Renyi}) holds. In particular, we study the $\gamma$-quasi clique and the logarithm of the maximum degree. For each of these properties, we verify (A\ref{item:property_continuous})--(A\ref{item:property_slowly_varying}) and deduce their asymptotic behavior for the IRG thanks to Theorem \ref{thm:asymptotics_IRG}. 

\subsection{The $\gamma$-quasi-clique}\label{sec:quasi-clique}
Given a graph $G=(V,E)$, let $G[S]$ be a sub-graph of $G$ induced by $S$ where $S \subset V$ is a vertex subset of $G$. We define a \emph{clique} as a subset of vertices $C \subset V$ such that $G[C]$ is a complete graph. For $\gamma \in [0,1]$, we define a $\gamma$-\emph{quasi-clique} as a subset of vertices $Q \subset V$ such that $G[Q]$ contains at least $\gamma \binom{|Q|}{2}$ edges. The $\gamma$-\emph{quasi-clique number} of $G$ is the size of the largest $\gamma$-\emph{quasi-clique} and we denote it with $\omega^\gamma(G)$. 
The $\gamma$-quasi-clique number satisfies (H\ref{item:f_increasing_partial_order}) since if we consider a graph $G$ and we add edges to $G$, then this property can only increase. Moreover, if we add zero-degree vertices to $G$, the $\gamma$-quasi-clique number remains the same, so it also satisfies (H\ref{item:f_decrasing_zero_degree}). 
Assumption (H\ref{item:asymptotic_property_Erdos_Renyi}) is satisfied thanks to \cite[Theorem 1]{balister2019dense}. Indeed, in \cite[Theorem 1]{balister2019dense} it is proven that for an ER random graph $\erdosrenyi{n,p}$ the following estimate holds for $\gamma<1$
\begin{align}
    \omega^\gamma(G) & >\frac{2}{D(\gamma,p)}(\log(n)-\log\log(n)+\log(eD(\gamma,p)/2) - \varepsilon  \\
    \omega^\gamma(G) & <\frac{2}{D(\gamma,p)}(\log(n)-\log\log(n)+\log(eD(\gamma,p)/2) + 1 + \varepsilon
\end{align}
where $D(\gamma,p)$ is the Kullback-Leibler divergence between the Bernoulli distributions \emph{Ber}$(\gamma)$ and \emph{Ber}$(p)$, i.e., 
\begin{equation}
    D(\gamma,p):= 
    \begin{cases}
    \gamma \log (\frac{\gamma}{p})+(1-\gamma)\log(\frac{1-\gamma}{1-p}), & \qquad \text{if } \gamma<1, \\
    \log(\frac{1}{p}), & \qquad \text{if } \gamma=1.
    \end{cases}
\end{equation}

Given $p^{\text{max}}=k(m,m)$ the maximal value of the kernel such that $p^{\text{max}}< \gamma \leq 1$ and $p^{\text{max}}_n= p^{\text{max}}/\lambda_n$,
it is easy to see that the function
\begin{equation}
    \omega_{n}^\gamma:=\frac{2\log(n)}{D(\gamma,p^{\text{max}}_n)}
\end{equation}
satisfies (A\ref{item:property_continuous}) and (A\ref{item:property_increasing}). So, we have to prove (A\ref{item:property_slowly_varying}) and to this end we take $r(n) := 1/\log(n)$. Then,
\begin{equation}
    \lim_{n \to \infty} \sup_{p\in(0,1)}\frac{\omega_{nr(n)}^\gamma}{\omega_n^\gamma} 
    %=\lim_{n \to \infty} \frac{2\log(nr(n))}{D(\gamma,p^{\text{max}})} \frac{D(\gamma,p^{\text{max}})}{2\log(n)}
    =1+\lim_{n \to \infty}\frac{\log(r(n))}{\log(n)}\to1.
\end{equation}
Thus, we have verified that all the conditions of Theorem \ref{thm:asymptotics_IRG} hold, and this concludes the proof of Theorem \ref{thm:cliques}.

\subsection{The logarithm of the maximum degree}\label{sec:Mdegree}
In this section, we verify that all the conditions for the \textit{logarithm} of the maximum degree of a dense ER random graph are satisfied. The \emph{maximum degree} $\Delta(G)$ of a graph $G$ is defined as $\Delta(G)=\max_{v\in V} d(v)$. The maximum degree satisfies (H\ref{item:f_increasing_partial_order}) since adding edges to a graph $G$ can only increase $\Delta(G)$. The maximum degree also satisfies (H\ref{item:f_decrasing_zero_degree}), since adding zero-degree vertices to a graph $G$ leaves $\Delta(G)$ unchanged. Since the maximum degree satisfies properties (H\ref{item:f_increasing_partial_order}), (H\ref{item:f_decrasing_zero_degree}), and since the logarithm is a monotone increasing function, the same properties are also satisfied by the logarithm of the maximum degree. 

Regarding (H\ref{item:asymptotic_property_Erdos_Renyi}), in \cite{bollobas1998random, durrett2010random, frieze2015introduction, frieze2024maximum} several estimates are given of the asymptotic behavior of the maximal degree.
In particular, in \cite[Theorem 3.5]{frieze2015introduction} it is proven that for an ER random graph $\erdosrenyi{n,p}$ with $p\in(0,1)$, for any fixed $\varepsilon>0$, with high probability
\begin{align*}
d_{n,p}(1-\varepsilon) < \Delta(G_{n,p}) < d_{n,p}(1+\varepsilon),
\end{align*}
where 
\begin{align}
    d_{n,p} :=(n-1)p.
\end{align} 
Thus, we note that (A\ref{item:property_continuous}) and (A\ref{item:property_increasing}) are trivially satisfied for $\log(d_n)$. To prove (A\ref{item:property_slowly_varying}) we take $r(n) := 1/\log(n)$, and compute
\begin{align}
    \sup_{p \in (0,1/\lambda_n)} \left |\frac{\log(d_{nr(n),p})}{\log(d_{n,p})} -1 \right |&=\sup_{p \in (0,1/\lambda_n)}  \frac{\log(n/\log(n)-1) +\log(p)}{\log(n-1) +\log(p)}-1  \\
    &= \frac{\log(n/\log(n)-1) -\log(\lambda_n)}{\log(n-1) -\log(\lambda_n)}-1 \to 0,
\end{align}
as $n\to\infty$, where we used that $p\mapsto\log(d_{nr(n),p})/\log(d_{n,p})$ is increasing.
The same conclusion holds when $p_n= (\omega_n \log(n))/n$ with $\omega_n \to \infty$ as $n \to \infty$ such that $p_n\in(0,1)$. Indeed, in this regime, the result follows from \cite[Theorem 3.2-(ii)]{frieze2015introduction}.

Since all assumptions of Theorem \ref{thm:asymptotics_IRG} are satisfied, we conclude that
\begin{align*}
\mathbb P(\log(d_{n,p}) (1-\varepsilon)< \log(\Delta(I_n(k_n))<\log(d_{n,p})(1+\varepsilon_n))\to1,
\end{align*}
as $n\to\infty$. In the following Theorem we phrase this result directly in terms of the maximal degree of the inhomogeneous random graph:
\begin{thm}[Maximum degree $\Delta$ of dense inhomogeneous random graphs]
Under the assumptions above, if $\lambda_n\equiv 1$ or $\lambda_n=n/( \omega_n \log(n))$ with $\omega_n \to \infty$ as $n \to \infty$ such that $\lambda_n\geq1$, then for any $\varepsilon>0$ we have 
\begin{align}%
\mathbb P \left (\Delta(\inhomogeneous{n}) \in [d_{n,p}^{1-\varepsilon}, d_{n,p}^{1 + \varepsilon}] \right ) \to 1,
\end{align}%
as $n\to\infty$, with $d_{n,p}=(n-1)p$.
\end{thm}
\section{Applications for sparse regime}\label{sec:applications_sparse}
In this section, we study the sparse regime, considering the maximum degree as a representative property $f$. Using its known asymptotic behavior for sparse ER random graphs, we verify (A\ref{item:property_continuous})–(A\ref{item:property_slowly_varying}) and apply Theorem \ref{thm:asymptotics_IRG}.
\subsection{The maximum degree}
In this section, we verify that all the required conditions are satisfied by the maximum degree of a sparse ER random graph. In Section \ref{sec:Mdegree}, we showed that the maximum degree satisfies properties (H\ref{item:f_increasing_partial_order}) and (H\ref{item:f_decrasing_zero_degree}).
%The \emph{maximum degree} $\Delta(G)$ of a graph $G$ is defined as $\Delta(G)=\max_{v\in V} d(v)$. The maximum degree satisfies (H\ref{item:f_increasing_partial_order}) since adding edges to a graph $G$ can only increase $\Delta(G)$. The maximum degree also satisfies (H\ref{item:f_decrasing_zero_degree}), since adding zero-degree vertices to a graph $G$ leaves $\Delta(G)$ unchanged. 
Regarding (H\ref{item:asymptotic_property_Erdos_Renyi}), in \cite[Theorem 3.4]{frieze2015introduction} it is proven that for an ER graph $\erdosrenyi{n,p}$ with $p=c/n$ and $c>0$, for any fixed $\varepsilon>0$, with high probability
\begin{align*}
\Delta_n(1-\varepsilon) < \Delta(G_{n,p}) < \Delta_n(1+\varepsilon),
\end{align*}
where, we recall, $\Delta_n := \log(n)/\log(\log(n))$. Note that this estimate does not depend on the constant $c$ and for this reason (A\ref{item:property_continuous}) and (A\ref{item:property_increasing}) are trivially satisfied.  To prove (A\ref{item:property_slowly_varying}) we take $r(n) := 1/\log(n)$, and
\begin{align}
    \frac{\Delta_{nr(n)}}{\Delta_n}=\frac{\log(nr(n))/\log(\log(nr(n)))}{\log(n)/\log(\log(n))}=\left ( 1-\frac{\log(\log(n))}{\log(n)} \right ) \left ( \frac{1}{1-\frac{\log(\log(\log(n)))}{\log(\log(n))}}\right ) \to 1,
\end{align}
as $n\to\infty$. Since all assumptions of Theorem \ref{thm:asymptotics_IRG} are satisfied, this concludes the proof of Theorem \ref{thm:degree}.

\section{Discussion}\label{sec:discussion}
\paragraph{A counterexample.}
Given a graph $G$, the distance $d(u,v)$ between two vertices $u$ and $v$ is the length of a shortest path joining $u$ and $v$, if it exists. If $G$ is a connected graph, then the average distance of $G$ is the average over all distances $d(u, v)$ for $u, v \in G$. If $G$ is not connected, then the average distance is customarily defined to be the average among all distances $d(u, v)$ where $u$ and $v$ belong to the same connected component.  
It turns out that the average distance is monotone only for connected graphs. Indeed, adding an edge that connects two different connected components may \textit{increase} the average distance, so  (H\ref{item:f_increasing_partial_order}) is not satisfied in general. 

Intuitively, our method works when the considered property is determined by the set of vertices with the most connections between them (i.e., those with weights close to the maximum of $k$). However, the average distance depends crucially on the entire topology of the graph, and for this reason, our method is not applicable in this context.

\paragraph{Maximum outside of the diagonal.} Our assumption that the maximum of $k$ lies on the diagonal is necessary in order to exploit the known results on the ER random graph. This is due to the way the lower bound is constructed. Recall that the vertex set of the corresponding graph $\tilde{G}_n$ is made up of only those vertices with weights in a small neighborhood of the maximum, see \eqref{eq:vertex_set_graph_lower_bound}. If $k(s,t)\leq k(n,m)$, for all $s,t\in[0,1]$ and $n\neq m$, then the vertices of the graph $\tilde{G}_n$ are $V' = \{i \in V~\vert~ W_i \in[m-r(n),m+r(n)] \cup [n-r(n),n+r(n)] \}$, and the edges are selected as in \eqref{coupling00}. However, the resulting graph $\tilde{G}_n = (V',E')$ is not distributed as an ER random graph. This is because the edges between two vertices with weight close to $m$ (resp.~close to $n$) have zero probability of being present. If the vertices are selected according to a different scheme (e.g., including the aforementioned missing edges with probability $p_{n}^{\text{inf}}$), then the resulting graph is not included in $G_n$, and thus it is not useful for computing a lower bound of $f(G_n)$. In other words, when the maximum lies outside of the diagonal, our technique relates asymptotics for the IRG with asymptotics for \textit{random bipartite graphs}.

\bibliographystyle{abbrv}
\bibliography{BJdoublecoupling}

\end{document}